\documentclass[11pt]{article}
\newcommand{\dated}{\mbox{} \hfill {\small [{\tt \today}]}} 
\newtheorem{theorem}{Theorem}[section]
\newtheorem{lemma}[theorem]{Lemma}

\newtheorem{proposition}[theorem]{Proposition}
\newtheorem{df}[theorem]{Definition}
\newenvironment{definition}{\begin{df} \rm}{\end{df}}

 \usepackage{amsmath,amssymb,amscd}
\newcommand{\pf}[1]{\trivlist \item[\hskip\labelsep\it #1\ ]}
\newcommand{\varpf}[1]{\trivlist \item[\hskip\labelsep\sc #1:]}
\newcommand{\qedbox}{$\rlap{$\sqcap$}\sqcup$}
\newcommand{\qed}{\qquad \qedbox \endtrivlist}
\newcommand{\varqed}{\hfill \rule{0.6em}{0.6em} \endtrivlist}
\newenvironment{proof}{\pf{Proof}}{\qed}

\newenvironment{remark}{\pf{Remark}}{\endtrivlist}
\newenvironment{remarks}{\pf{Remarks} 
   \begin{enumerate}}{\end{enumerate} \endtrivlist}
\newenvironment{example}{\pf{Example}}{\endtrivlist}
\newenvironment{examples}{\pf{Examples} 
   \begin{enumerate}}{\end{enumerate} \endtrivlist}
\newenvironment{items}{
  \begin{enumerate} 
                    
  }{\end{enumerate}}

\newenvironment{keywords}{\noindent\small {\it Keywords\/}:}{\vskip 4pt}
\newenvironment{classification}{\noindent\small 2000 {\it Mathematics Subject
Classification\/}:}{\vskip 12pt}

\renewcommand{\iff}{\quad\Longleftrightarrow\quad}

\newcommand{\comps}{{\mathbb C}}
\newcommand{\reals}{{\mathbb R}}

\newcommand{\free}{{\mathbb F}}

\newcommand{\void}{\varnothing}
\newcommand{\tensor}{\otimes}

\newcommand{\Tensor}{\hat{\otimes}}

\newcommand{\cstar}{{C^\ast}}
\newcommand{\wstar}{{W^\ast}}

\newcommand{\id}{{\mathrm{id}}}

\newcommand{\A}{{\mathfrak A}}

\newcommand{\M}{{\mathfrak M}}

\newcommand{\VN}{\operatorname{VN}}

\newcommand{\AP}{\operatorname{AP}}
\newcommand{\WAP}{\operatorname{WAP}}
\newcommand{\PIF}{{\operatorname{\cal PIF}}}
\title{Operator amenability of Fourier--Stieltjes algebras}
\author{{\it Volker Runde}\thanks{Research supported by NSERC under grant no.\ 227043-00.} \and {\it Nico Spronk}}
\date{}
\begin{document}
\maketitle
\begin{abstract}
In this paper, we investigate, for a locally compact group $G$, the operator amenability of the Fourier-Stieltjes algebra $B(G)$ and of the reduced Fourier-Stieltjes algebra $B_r(G)$. The natural
conjecture is that any of these algebras is operator amenable if and only if $G$ is compact. We partially prove this conjecture with mere operator amenability replaced by operator $C$-amenability 
for some constant $C < 5$. In the process, we obtain a new decomposition of $B(G)$, which can be interpreted as the non-commutative counterpart of the decomposition of $M(G)$ into the discrete and 
the continuous measures. We further introduce a variant of operator amenability --- called operator Connes-amenability --- which also takes the dual space structure on $B(G)$ and $B_r(G)$ into account. 
We show that $B_r(G)$ is operator Connes-amenable if and only if G is amenable. Surprisingly, $B(\free_2)$ is operator Connes-amenable although $\free_2$, the free group in two generators, fails to be amenable.
\end{abstract}
\section*{Introduction}
In his now classic memoir \cite{Joh1}, B.\ E.\ Johnson initiated the theory of amenable Banach algebras. The choice of terminology is motivated by
\cite[Theorem 2.5]{Joh1}: a locally compact group $G$ is amenable if and only if its group algebra $L^1(G)$ is an amenable Banach algebra. There are other 
Banach algebras associated with a locally compact $G$ which are as natural objects of study as $L^1(G)$, e.g.\ the measure algebra $M(G)$. If $G$ is discrete
and amenable, then $M(G) = \ell^1(G) = L^1(G)$ is amenable by Johnson's theorem. It was conjectured by A.\ T.-M.\ Lau and R.\ J.\ Loy that $M(G)$ is amenable
{\it only if\/} $G$ is discrete and amenable (\cite{LL}), a conjecture that was ultimately confirmed by H.\ G.\ Dales, F.\ Ghahramani, and A.\ Ya.\ Helemski\u{\i}
(\cite{DGH}).
\par
In \cite{Eym}, P.\ Eymard introduced, for an arbitrary locally compact $G$, its Fourier algebra $A(G)$ and its Fourier--Stieltjes algebra $B(G)$. If $G$ is abelian with dual group
$\hat{G}$, then the Fourier and Fourier--Stieltjes transform, respectively, yield $A(G) \cong L^1(\hat{G})$ and $B(G) \cong M(\hat{G})$. Disappointingly, the amenability of $A(G)$ reflects the
amenability of $G$ rather inadequately: there are compact groups $G$, e.g.\ $G = \operatorname{SO}(3)$, for which $A(G)$ fails to be amenable (\cite{Joh3}). It would seem that the only locally
compact groups $G$ for which $A(G)$ is known to be amenable are those which have a closed, abelian subgroup with finite index (\cite{LLW} or \cite{For}).
\par
Being the predual of the group von Neumann algebra $\VN(G)$, the Fourier algebra $A(G)$ has a canonical operator space structure. In \cite{Rua}, Z.-J.\ Ruan introduced a variant of amenability
--- called {\it operator amenability\/} --- which takes the operator space structure of $A(G)$ into account. As it turns out, operator amenability is the ``right'' notion of amenability
for $A(G)$ in the sense that it characterizes the amenable, locally compact groups: $A(G)$ is operator amenable if and only if $G$ is amenable (\cite[Theorem 3.6]{Rua}).
\par
Let $\cstar(G)$ and $C^\ast_r(G)$ denote the full and the reduced group $\cstar$-algebra of $G$, respectively. Then $B(G) = C^\ast(G)^\ast$ and the reduced Fourier--Stieltjes algebra
$B_r(G) = C^\ast_r(G)^\ast$ also have canonical operator space structures turning them into completely contractive Banach algebras. It is thus natural to ask for which $G$, the algebras $
B(G)$ and $B_r(G)$, respectively, are operator amenable (\cite[Problem 32]{LoA}). Since $B_r(G) = B(G) = M(\hat{G})$ for abelian $G$, \cite[Theorem 1.1]{DGH} suggests that this is the case 
if and only if $G$ is compact. We have not been able to prove this conjecture in full. However, if we replace mere operator amenability by what we shall call operator $C$-amenability: 
the Fourier--Stieltjes algebra $B(G)$ --- and, equivalently, $B_r(G)$ --- is operator $C$-amenable for some $C < 5$ if and only if $G$ is compact.
\par
In \cite{LoA}, it was conjectured that, if we want to capture the amenability of a locally compact group $G$ in terms of an amenability condition for $B(G)$ or $B_r(G)$, this notion of
amenability needs to take both the operator space and the dual space structure of $B(G)$ and $B_r(G)$ into account. We introduce such a notion --- called {\it operator Connes-amenability\/} --- and show that, indeed,
$B_r(G)$ is operator Connes-amenable if and only if $G$ is amenable. Surprisingly, there are non-amenable, locally compact groups $G$ --- including $\free_2$ --- for which $B(G)$ is operator Connes-amenable.
\section{Completely contractive Banach algebras and operator amenability}
Since there are now several expository sources on the theory of operator spaces available (\cite{ER}, \cite{Pis}, and \cite{Wit}), we refrain from introducing the basics of operator space theory.
We will adopt the notation from \cite{ER}; in particular, $\Tensor$ stands for the projective tensor product of operator spaces and not of Banach spaces.
\par
We briefly recall a few definitions and results from \cite{Rua}.
\begin{definition}
A Banach algebra $\A$ which is also an operator space is called {\it completely contractive\/} if the multiplication of $\A$ is a completely contractive bilinear map.
\end{definition}
\par
Clearly, $\A$ is completely contractive if and only if the multiplication of $\A$ induces a complete contraction $\Delta \!: \A \Tensor \A \to \A$.
\begin{examples}
\item For any Banach algebra $\A$, the maximal operator space $\max \A$ is completely contractive.
\item If $\mathfrak H$ is a Hilbert space, then any closed subalgebra of ${\cal B}({\mathfrak H})$ is completely contractive.
\item We denote the $W^\ast$-tensor product by $\bar{\tensor}$. A {\it Hopf--von Neumann algebra\/} is a pair $(\M, \nabla)$, where $\M$ is a von Neumann algebra, and 
$\nabla$ is a {\it co-multiplication\/}: a unital, $w^\ast$-continuous, and injective $^\ast$-homomorphism $\nabla \!: \M \to \M \bar{\tensor} \M$ which is co-associative, i.e.\ the diagram
\[
  \begin{CD}
  \M  @>\nabla>> \M \bar{\tensor} \M \\
  @V{\nabla}VV                                          @VV{\nabla \tensor \id_\M}V                     \\
  \M \bar{\tensor} \M @>>{\id_\M \tensor \nabla}> \M \bar{\tensor} \M \bar{\tensor} \M
  \end{CD}    
\]
commutes. Let $\M_\ast$ denote the unique predual of $\M$. By \cite[Theorem 7.2.4]{ER}, we have $\M \bar{\tensor} \M \cong (\M_\ast \Tensor \M_\ast)^\ast$. Thus $\nabla$ induces a
complete contraction $\nabla_\ast \!: \M_\ast \Tensor \M_\ast \to \M_\ast$ turning $\M_\ast$ into a completely contractive Banach algebra. 
\item Let $G$ be a locally compact group, and let $\wstar(G) := C^\ast(G)^{\ast\ast}$. There is a canonical $w^\ast$-continuous unitary representation $\omega \!: G \to \wstar(G)$, the {\it universal representation\/} of $G$,
with the following universal property: for any representation (always WOT-continuous and unitary) $\pi$ of $G$ on a Hilbert space, there is unique $w^\ast$-continuous $^\ast$-homomorphism 
$\theta \!: \wstar(G) \to \pi(G)''$ such that $\pi = \theta \circ \omega$. Applying this universal property to the representation
\[
  G \to \wstar(G) \bar{\tensor} \wstar(G), \quad x \mapsto \omega(x) \tensor \omega(x)
\]
yields a co-multiplication $\nabla \!: \wstar(G) \to \wstar(G) \bar{\tensor} \wstar(G)$. Hence, $B(G) := \cstar(G)^\ast$, the {\it Fourier--Stieltjes\/} algebra of $G$, is a completely contractive Banach algebra. 
Since $B_r(G)$ and $A(G)$ are closed ideals of $B(G)$ (see \cite{Eym}), they are also completely contractive Banach algebras. (It is not hard to see that the operator space structures on $B_r(G)$ and $A(G)$ inherited from 
$B(G)$ coincide with those they have as the preduals of $C^\ast_r(G)^{\ast\ast}$ and $\VN(G)$, respectively.)
\end{examples}
\begin{definition}
Let $\A$ be a completely contractive Banach algebra. An operator $\A$-bimodule $E$ is an $\A$-bimodule $E$ which is also an operator space such that the module actions 
\[
  \A \times E \to E, \quad (a,x) \mapsto a \cdot x 
  \qquad\text{and}\qquad
  E \times \A \to E, \quad (x, a) \mapsto x \cdot a
\]
are completely bounded.
\end{definition}
\par
Similarly, one defines left and right operator $\A$-modules. If $E$ is a left and $F$ is a right operator $\A$-module, then $E \Tensor F$ becomes an operator $\A$-bimodule in a canonical fashion
via
\[
  a \cdot (x \tensor y) := a \cdot x \tensor y \quad\text{and}\quad\text (x \tensor y)\cdot a := x \tensor y \cdot a \qquad (a \in \A, \, x \in E, \, y \in F).
\]
In particular, $\A \Tensor \A$ is an operator $\A$-bimodule in a canonical way.
\par
For any operator $\A$-bimodule $E$, its dual module $E^\ast$ is also an operator $\A$-bimodule. We shall call an operator $\A$-bimodule $E$ {\it dual\/} if it is of the form
$E = (E_\ast)^\ast$ for some operator $\A$-bimodule $E_\ast$.
\begin{definition}
A completely contractive Banach algebra $\A$ is called {\it operator amenable\/} if every completely bounded derivation from $\A$ into a dual operator $\A$-bimodule is inner.
\end{definition}
\par
There is an intrinsic characterization of amenable Banach algebras in terms of approximate diagonals (\cite{Joh2}). This characterization has an analogue for operator amenable, completely
contractive Banach algebras (\cite[Proposition 2.4]{Rua}):   
\begin{theorem} \label{diagthm}
The following are equivalent for a completely contractive Banach algebra $\A$:
\begin{items}
\item $\A$ is operator amenable.
\item There is an\/ {\rm approximate operator diagonal} for $\A$, i.e.\ a bounded net $( {\mathbf d}_\alpha )_\alpha$ in $\A \Tensor \A$ such that
\[
  a \cdot {\mathbf d}_\alpha - {\mathbf d}_\alpha \cdot a \to 0 \quad\text{and}\quad a \Delta {\mathbf d}_\alpha \to a \qquad (a \in \A).
\]
\item There is a\/ {\rm virtual operator diagonal} for $\A$, i.e.\ an element ${\mathbf D} \in (\A \Tensor \A)^{\ast\ast}$ such that
\[
  a \cdot {\mathbf D} = {\mathbf D} \cdot a \quad\text{and}\quad a \Delta^{\ast\ast} {\mathbf D} = a \qquad (a \in \A).
\]
\end{items}
\end{theorem}
\par
In analogy with the classical situation, Theorem \ref{diagthm} allows for a refinement of the notion of operator amenability:
\begin{definition}
Let $C \geq 1$. A completely contractive Banach algebra $\A$ is called {\it operator $C$-amenable\/} if there is an approximate operator diagonal for $\A$ bounded by $C$.
\end{definition}
\begin{example}
For any amenable, locally compact group $G$, the Fourier algebra $A(G)$ is operator $1$-amenable (this is implicitly shown in \cite{Rua}).
\end{example}
\par
We conclude this preliminary section with a lemma, which is the operator analogue of a classical result (\cite[Theorem 2.3.7]{LoA}); given Theorem \ref{diagthm}, the proof from \cite{LoA} carries over with the 
obvious modifications:
\begin{lemma} \label{ideal}
Let $\A$ be an operator amenable, completely contractive Banach algebra. Then the following are equivalent for a closed ideal $I$ of $\A$:
\begin{items}
\item $I$ is operator amenable.
\item $I$ has a bounded approximate identity.
\item $I$ is completely weakly complemented, i.e.\ there is a completely bounded projection from $\A^\ast$ onto $I^\perp$.
\end{items}
\end{lemma}
\begin{remark}
Of course, (iii) is satisfied whenever $I$ is completely complemented, i.e.\ if there is a completely bounded projection from $\A$ onto $I$.
\end{remark}
\section{A decomposition for $B(G)$}
Let $G$ be a locally compact group. Then we have a direct sum decomposition $M(G) = \ell^1(G) \oplus M_c(G)$, where $M_c(G)$ denotes the ideal of continuous measures in $M(G)$. This
decomposition was crucial in the proof of \cite[Theorem 1.1]{DGH}. In this section, we establish an analogous decomposition for $B(G)$.
\par
Let $G$ be a abelian with dual group $\hat{G}$ whose Bohr compactification we denote by $b\hat{G}$; we write $G_d$ for the group $G$ equipped with the discrete topology. For $\mu \in M(G)$, we denote its Fourier--Stieltjes 
transform in $B(\hat{G})$ by $\hat{\mu}$. Then we have for $\mu \in M(G)$:
\begin{eqnarray*}
  \mu \in \ell^1(G) & \iff & \mu \in M(G_d) \\
  & \iff & \hat{\mu} \in B( \widehat{G_d} ) \\
  & \iff & \hat{\mu} \in B( b \hat{G} ) \\
  & \iff & \text{$\hat{\mu} \in B(\hat{G})$ is almost periodic},
\end{eqnarray*}
where the last equivalence holds by \cite[(2.27) Corollaire 4]{Eym}.
\par
This suggests that the appropriate replacement for $\ell^1(G)$ in the Fourier--Stieltjes algebra context is $B(G) \cap \AP(G)$, where $\AP(G)$ denotes the algebra of all almost periodic functions on $G$.
It is well known (see \cite[3.2.16]{Pal}, for example) that $\AP(G)$ is a commutative $\cstar$-algebra whose character space is a compact group denoted by $aG$ (for abelian $G$, we have $aG = bG$). We will 
first give an alternative description of $B(G) \cap \AP(G)$ which will turn out to be useful later on.
\par
Let $G$ be a locally compact group, and let $\cal R$ be any family of representations of $G$. We denote by $A_{\cal R}(G)$ the closed linear span in $B(G)$ of the coefficient functions of all
representations in $\cal R$, i.e.\ of all functions of the form
\[
  G \to \comps, \quad x \mapsto \langle \rho(x) \xi, \eta \rangle,
\]
where $\rho \in \cal R$, and $\xi$ and $\eta$ are vectors in the corresponding Hilbert space. If $\cal R$ is the family of all representations of $G$, then $A_{\cal R}(G) = B(G)$, and if $\cal R$ just consists of
the left regular representation, then $A_{\cal R}(G)= A(G)$. Let $\cal F$ denote the family of all finite-dimensional representations of $G$. Since $\cal F$ is closed under
taking tensor products, it is immediate that $A_{\cal F}(G)$ is a (completely contractive) Banach algebra.
\begin{proposition} \label{nicoprop}
Let $G$ be a locally compact group. Then $A_{\cal F}(G) = B(G) \cap \AP(G)$, and we have a canonical completely isometric isomorphism between $A_{\cal F}(G)$ and $B(aG)$.
\end{proposition}
\begin{proof}
In view of \cite[(2.27) Corollaire 4]{Eym}, it is sufficient to prove the second assertion only.
\par
Let $\iota \!: G \to a G$ denote the (not necessarily injective) canonical map. It is easy to
see that $A_{\cal F}(G) \cong B(aG)$ via
\begin{equation} \label{AP}
  B(aG) \to B(G), \quad f \mapsto f \circ \iota.
\end{equation}
We claim that (\ref{AP}) is a complete isometry. To see this, let $\omega_G \!: G \to \wstar(G)$ and $\omega_{aG} \!: aG \to \wstar(aG)$ denote the universal representations of $G$ and $aG$, respectively. Applying the
universal property of $\omega_G \!: G \to \wstar(G)$ to $\omega_{aG} \circ \iota \!: G \to \wstar(aG)$ yields a (necessarily surjective) $w^\ast$-continuous $^\ast$-homomorphism $\pi \!: \wstar(G) \to \wstar(aG)$. It
is immediate that (\ref{AP}) is the adjoint of $\pi$. Hence, (\ref{AP}) is a complete isometry by \cite[Theorem 4.1.8]{ER}.
\end{proof}
\begin{remarks}
\item Note that $A_{\cal F}(G)$ can be very small relative to $B(G)$: for example, if $G = \operatorname{SL}(2,\reals)$, we have, $A_{\cal F}(G) = \comps$.
\item Suppose that $G$ is non-compact. Since $B(G)$ is a complete invariant for $G$ (\cite[Corollary]{Wal}), it follows that $B(G) \not\cong B(aG)$ and thus $A_{\cal F}(G) \subsetneq B(G)$ by Proposition \ref{nicoprop}.
\end{remarks}
\par
Let $G$ be a locally compact group. For any function $f$ on $G$ and $x \in G$, we define the left and the right translate of $f$ by $x$ by letting
\[
  (L_x f)(y) := f(xy) \quad\text{and}\quad (R_x f)(y) := f(yx) \qquad (y \in G).
\]
A linear space $E$ of functions on $G$ is said to be {\it translation invariant\/} if $L_x f, R_x f \in E$ for all $f \in E$ and $x \in G$.
\par
We record the following well known lemma for convenience:
\begin{lemma} \label{translem1}
Let $G$ be a locally compact group. Then the following are equivalent for a closed subspace $E$ of $B(G)$:
\begin{items}
\item $E$ is translation invariant;
\item $E$ is a $\wstar(G)$-submodule of $B(G)$;
\item $E = p \cdot B(G)$ for a unique central projection $p \in \wstar(G)$.
\end{items}
\end{lemma}
\begin{proof}
(i) $\Longleftrightarrow$ (ii) is \cite[Proposition 1.(i)]{Wal}, and (ii) $\Longleftrightarrow$ (iii) is a well known general fact about von Neumann algebras (which can be found in \cite{Tak}, for instance).
\end{proof}
\par
Let $\cal R$ be a family of representations of $G$. Then it is clear that $A_{\cal R}(G)$ is translation invariant. Hence, there is a unique central projection $p_{\cal R} \in \wstar(G)$ such that 
$A_{\cal R}(G) = p_{\cal R} \cdot B(G)$. 
\par
For any representation $\pi$ of $G$, we denote its canonical $w^\ast$-continuous extension to $\wstar(G)$ by $\pi$ as well. We call a representation $\pi$ of $G$ {\it purely infinite-dimensional\/} if
$\pi(p_{\cal F}) = 0$. We denote the family of all purely infinite-dimensional representations of $G$ by $\PIF$; note that $\PIF \neq \void$ if $G$ is not compact.
\begin{theorem} \label{decompthm}
Let $G$ be a locally compact group. Then the following are equivalent and true:
\begin{items}
\item $A_\PIF(G)$ is an ideal of $B(G)$.
\item The map
\[
   B(G) \to A_{\cal F}(G), \quad f \mapsto p_{\cal F} \cdot f
\]
is an algebra homomorphism.
\item $\nabla p_{\cal F} = p_{\cal F} \tensor p_{\cal F}$.
\end{items}
\end{theorem}
\begin{proof}
It is immediately checked that (i) and (ii) are equivalent.
\par
Let $x \in \wstar(G)$ and $f,g \in B(G)$, and note that
\begin{eqnarray*}
  \langle x, p_{\cal F} \cdot (fg) - (p_{\cal F} \cdot f) (p_{\cal F} \cdot g) \rangle & = & \langle x p_{\cal F}, fg \rangle - \langle x, (p_{\cal F} \cdot f) (p_{\cal F} \cdot g) \rangle \\
  & = & \langle \nabla(x p_{\cal F}), f \tensor g \rangle - \langle \nabla x, (p_{\cal F} \cdot f) \tensor (p_{\cal F} \cdot g) \rangle \\ 
  & = & \langle \nabla(x) \nabla(p_{\cal F}) - (\nabla x) (p_{\cal F} \tensor p_{\cal F}), f \tensor g \rangle \\
  & = & \langle \nabla(x) (\nabla p_{\cal F} - p_{\cal F} \tensor p_{\cal F}), f \tensor g \rangle.
\end{eqnarray*}
This proves the equivalence of (ii) and (iii).
\par
We shall now verify that (iii) is indeed true. 
\par
Let $\WAP(G)$ denote the weakly almost periodic functions on $G$ (see \cite{Bur} for the definition of $\WAP(G)$ and further information). By \cite[Theorem 3.1]{Bur}, we have $B(G) \subset \WAP(G)$. Taking the adjoint of this
inclusion map, we obtain a canonical map $\pi \!: \WAP(G)^\ast \to \wstar(G)$. Since $\WAP(G)$ is an introverted subspace of $\ell^\infty(G)$, its dual $\WAP(G)^\ast$ is an Banach algebra in a canonical manner. It is routinely
verified --- e.g.\ by checking multiplicativity on $M(G)$ --- that $\pi$ is a $^\ast$-homomorphism. The character space $wG$ of $\WAP(G)$ is a compact, semitopological semigroup containing a topologically isomorphic copy of $G$. 
The kernel $K(wG)$ of $wG$ is intersection of all ideals of $wS$; it is non-empty by \cite[Theorems 2.1 and 2.2]{Bur}, and by \cite[Theorems 2.7]{Bur}, it is a compact group.  Let $e_{K(wG)}$ denote its identity element. Then
by (the proof of) \cite[Theorem 2.22]{Bur}, we have $\AP(G) = e_{K(wG)} \cdot \WAP(G)$. It follows that $p_{\cal F} = \pi(e_{K(wG)})$. In particular, $p_{\cal F}$ is a character on $B(G)$. By \cite[Theorem 1.(ii)]{Wal}, this
implies (iii).
\end{proof}
\begin{remarks}
\item Let $G$ be a non-discrete locally compact group. Then we have a further decomposition of $M_c(G)$, namely $M_c(G) = M_s(G) \oplus  L^1(G)$, where $M_s(G)$ denotes the measures in $M_c(G)$ which are singular with respect to
left Haar measure. The decomposition of $M(G)$ into $\ell^1(G) \oplus M_s(G)$ and $L^1(G)$ has long been known to have a $B(G)$-analogue (see \cite{Ars} and \cite{Mia}). In view of Theorem \ref{decompthm}, we now have a complete
analogue for $B(G)$ of the decomposition of the measure algebra into its discrete part, its singular, continuous part, and its absolutely continuous part.
\item Let $G$ be a non-compact, locally compact group. Then $A(G)$ is a translation invariant subspace of $B(G)$ having zero intersection with $A_{\cal F}(G)$. It follows that $A(G) \subset A_{\PIF}(G)$. Since for
a non-discrete, locally compact group, the absolutely continuous measures are properly contained in the continuous measures, the natural conjecture is that  $A(G) \subsetneq A_{\PIF}(G)$. This conjecture seems to be open
for general locally compact groups, even in the amenable case.
\end{remarks}
\section{Operator non-amenability for $B(G)$ and $B_r(G)$ if $G$ is not compact}
We will now use Theorem \ref{decompthm} to show that $B(G)$ --- and, equivalently, $B_r(G)$ --- cannot be operator $C$-amenable with $C < 5$ unless $G$ is compact.
\par
We first need a purely operator space theoretic lemma.
\par
Given two opertor spaces $E_1$ and $E_2$, their operator space $\ell^\infty$-direct sum $E_1 \oplus_\infty E_2$ is defined by taking the Banach space $\ell^\infty$-direct sum on each matrix level. It is then
immediate that $E_1 \oplus_\infty E_2$ is again an operator space. If $F$ is another operator space, then it is immediately checked that
\begin{equation} \label{ellinfty}
  {\cal CB}(F, E_1 \oplus_\infty E_2) \cong {\cal CB}(F,E_1) \oplus_\infty {\cal CB}(F,E_2)
\end{equation}
canonically as Banach spaces. From the definition of the operator space structures on ${\cal CB}(E_1 \oplus_\infty E_2,F)$, ${\cal CB}(E_1,F)$, and ${\cal CB}(E_2,F)$ (see \cite[p.\ 45]{ER}), it follows
that the identification (\ref{ellinfty}) is even a complete isometry.
\par
The canonical embedding of $E_1 \oplus E_2$ into $(E_1^\ast \oplus_\infty E_2^\ast)^\ast$ equips $E_1 \oplus E_2$ with another operator space structure, denoted by $E_1 \oplus_1 E_2$. On the Banach space level,
this is just the ordinary $\ell^1$-direct sum of Banach spaces. Replacing $E_1^\ast$ and $E_2^\ast$ with $E_1$ and $E_2$, respectively, in (\ref{ellinfty}) and combining the duality result \cite[Corollary 7.1.5]{ER}
with the commutativity of $\Tensor$, we obtain:
\begin{lemma} \label{ellone}
Let $E_1$, $E_2$, and $F$ be operator spaces. Then we have a canonical completely isometric isomorphism
\[
  (E_1 \oplus_1 E_2) \Tensor F \cong  (E_1 \Tensor F) \oplus_1 (E_2 \Tensor F).
\]
\end{lemma}
\par
We can now prove the main result of this section:
\begin{theorem} \label{opamthm}
For a locally compact group, the following are equivalent:
\begin{items}
\item $G$ is compact.
\item $B_r(G)$ is operator $C$-amenable for some $C < 5$.
\item $B(G)$ is operator $C$-amenable for some $C < 5$.
\end{items}
\end{theorem}
\begin{proof}
(i) $\Longrightarrow$ (ii): If $G$ is compact, then $B_r(G) = B(G) = A(G)$. Since $A(G)$ is operator $1$-amenable, this proves (ii).
\par
(ii) $\Longrightarrow$ (iii): Since $A(G)$ is a closed $C^\ast_r(G)$-submodule of $B_r(G)$, there is a projection $p \in C^\ast_r(G)^{\ast\ast}$ such that $A(G) = p \cdot B_r(G)$. In particular, $A(G)$ is a completely complemented
ideal of $B_r(G)$ and thus operator amenable by Lemma \ref{ideal}. By \cite[Theorem 3.6]{Rua}, this implies the amenability of $G$ and thus $B_r(G) = B(G)$ by \cite[(4.21) Theorem]{Pat}. 
\par 
(iii) $\Longrightarrow$ (i): Assume towards a contradiction that $G$ is not compact. Let $( {\mathbf d}_\alpha )_{\alpha \in {\mathbb A}}$ be an approximate operator diagonal for $B(G)$ bounded by $C < 5$. Without loss of generality, 
suppose that $\Delta {\mathbf d}_\alpha = 1$ for all $\alpha \in \mathbb A$. We then have
\begin{equation} \label{diagdecomp}
  {\mathbf d}_\alpha = p_{\cal F} \cdot {\mathbf d}_\alpha \cdot p_{\cal F} +  p_{\cal F} \cdot {\mathbf d}_\alpha \cdot p_{\PIF} + p_{\PIF} \cdot {\mathbf d}_\alpha \cdot p_{\cal F}
  + p_{\PIF} \cdot {\mathbf d}_\alpha \cdot p_{\PIF} \qquad (\alpha \in {\mathbb A}).
\end{equation}
Since $B(G) = p_{\cal F} \cdot B(G) \oplus_1 p_{\PIF} \cdot B(G)$ in the operator space sense, Lemma \ref{ellone} and (\ref{diagdecomp}) yield
\begin{eqnarray} 
  \lefteqn{\| p_{\cal F} \cdot {\mathbf d}_\alpha \cdot p_{\cal F} \| + \| p_{\cal F} \cdot {\mathbf d}_\alpha \cdot p_{\PIF} \|} & & \nonumber \\ 
  & & \mbox{}  + \| p_{\PIF} \cdot {\mathbf d}_\alpha \cdot p_{\cal F} \|
  + \| p_{\PIF} \cdot {\mathbf d}_\alpha \cdot p_{\PIF} \| = \| {\mathbf d}_\alpha \| \leq C < 5 \qquad (\alpha \in {\mathbb A}). \label{diagdecomp2}
\end{eqnarray}
\par
First note that, by Theorem \ref{decompthm}, we have
\begin{equation} \label{Fpart}
  \Delta( p_{\cal F} \cdot {\mathbf d}_\alpha \cdot p_{\cal F} ) = p_{\cal F} \cdot \Delta {\mathbf d}_\alpha = p_{\cal F} \cdot 1 = 1.
\end{equation}
Since $\Delta$ is a (complete) contraction, this yields in turn that
\begin{equation} \label{estim1}
  \| p_{\cal F} \cdot {\mathbf d}_\alpha \cdot p_{\cal F} \| \geq 1 \qquad (\alpha \in {\mathbb A}).
\end{equation}
Let $\cal U$ be an ultrafilter on $\mathbb A$ that dominates the order filter. The we have for $f \in A(G) \subset A_{\PIF}(G)$:
\begin{eqnarray*}
  f (\text{$w^\ast$-}\lim_{\cal U} \Delta( {\mathbf d}_\alpha \cdot p_{\cal F} )) & = & \text{$w^\ast$-}\lim_{\cal U} \Delta( (f \cdot {\mathbf d}_\alpha) \cdot p_{\cal F} ) \\
  & = & \text{$w^\ast$-}\lim_{\cal U} \Delta( ({\mathbf d}_\alpha \cdot f) \cdot p_{\cal F} ) \\
  & = & 0.
\end{eqnarray*}
It follows that $\text{$w^\ast$-}\lim_{\cal U} \Delta( {\mathbf d}_\alpha \cdot p_{\cal F} ) = 0$. Combining this with (\ref{Fpart}), we obtain 
\[
  \text{$w^\ast$-}\lim_{\cal U} \Delta( p_{\PIF} \cdot {\mathbf d}_\alpha \cdot p_{\cal F} ) = -1
\]
and therefore
\begin{equation} \label{estim2}
  \lim_{\cal U} \|  p_{\PIF} \cdot {\mathbf d}_\alpha \cdot p_{\cal F} \| \geq \lim_{\cal U} \| \Delta( p_{\PIF} \cdot {\mathbf d}_\alpha \cdot p_{\cal F} ) \| \geq 1.
\end{equation}
Analoguously, we see that $\text{$w^\ast$-}\lim_{\cal U} \Delta( p_{\cal F} \cdot {\mathbf d}_\alpha \cdot p_{\PIF} ) = -1$ and consequently
\begin{equation} \label{estim3}
  \lim_{\cal U} \|  p_{\cal F} \cdot {\mathbf d}_\alpha \cdot p_{\PIF} \| \geq 1.
\end{equation}
Since
\begin{eqnarray*}
  1 & = & \text{$w^\ast$-}\lim_{\cal U}\Delta(p_{\cal F} \cdot {\mathbf d}_\alpha \cdot p_{\cal F}) \\
  & = & \text{$w^\ast$-}\lim_{\cal U} (\Delta( {\mathbf d}_\alpha) - \Delta(p_{\PIF} \cdot {\mathbf d}_\alpha \cdot p_{\cal F}) - \Delta(p_{\cal F} \cdot {\mathbf d}_\alpha \cdot p_{\PIF}) -
  \Delta(p_{\PIF} \cdot {\mathbf d}_\alpha \cdot p_{\PIF})) \\
  & = & 1 - \text{$w^\ast$-}\lim_{\cal U} \Delta( p_{\PIF} \cdot {\mathbf d}_\alpha \cdot p_{\cal F} ) - \text{$w^\ast$-}\lim_{\cal U} \Delta( p_{\cal F} \cdot {\mathbf d}_\alpha \cdot p_{\PIF} ) \\
  & & \mbox{} - \text{$w^\ast$-}\lim_{\cal U} \Delta( p_{\PIF} \cdot {\mathbf d}_\alpha \cdot p_{\PIF} ) \\
  & = & 3 - \text{$w^\ast$-}\lim_{\cal U} \Delta( p_{\PIF} \cdot {\mathbf d}_\alpha \cdot p_{\PIF} ),
\end{eqnarray*}
it follows that $\text{$w^\ast$-}\lim_{\cal U} \Delta( p_{\PIF} \cdot {\mathbf d}_\alpha \cdot p_{\PIF} ) = 2$. We thus obtain
\begin{equation} \label{estim4}
  \lim_{\cal U} \|  p_{\PIF} \cdot {\mathbf d}_\alpha \cdot p_{\PIF} \| \geq 2.
\end{equation}
\par
Altogether, (\ref{estim1}), (\ref{estim2}), (\ref{estim3}), and (\ref{estim4}) contradict (\ref{diagdecomp2}).
\end{proof}
\begin{remarks}
\item The proof of (ii) $\Longrightarrow$ (i) shows that, whenever $B_r(G)$ --- or, equivalently, $B(G)$ --- is operator amenable, then $G$ is amenable.
\item We strongly suspect that $B(G)$ and $B_r(G)$ are operator amenable only if $G$ is compact. One possible way of proving this would be to follow the route outlined (for measure algebras) in \cite{DGH}: assume that
$B(G)$ is operator amenable, but that $G$ is not compact. Then Lemma \ref{ideal} implies that $A_{\PIF}(G)$ is operator amenable and thus has a bounded approximate identity. This, in turn, would imply that every 
element of $A_{\PIF}(G)$ is a product of two elements in $A_{\PIF}(G)$ by Cohen's factorization theorem (\cite[5.2.4 Corollary]{Pal}). We believe that this is not true, but have been unable to confirm this belief with a proof.
\item Another open question related to Theorem \ref{opamthm} is for which locally compact groups $G$, the Fourier--Stieltjes algebra $B(G)$ is amenable in the classical sense. The corresponding question for the Fourier algebra 
is also still open: as mentioned in the introduction, the only locally compact groups $G$ for which $A(G)$ is known to be amenable are those with an abelian subgroup of finite index, and it is plausible to conjecture that these 
are indeed the only ones. The plausible conjecture for $B(G)$ is that it is amenable if and only if $G$ is compact and has an abelian subgroup of finite index.
\end{remarks}
\section{Operator Connes-amenability}
Amenability in the sense of \cite{Joh1} is not the ``right'' notion of amenability for von Neumann algebras because it is too restrictive to allow for the
development of a reasonably rich theory (\cite{Was}). In \cite{JKR}, a variant of amenability --- christened Connes-amenability in \cite{Hel} --- was introduced for
von Neumann algebras, which takes the normal structure in von Neumann algebras into account. This notion of amenability has turned out to be equivalent to a number
of important $W^\ast$-algebraic properties, such as injectivity and semidiscreteness; see \cite[Chapter 6]{LoA} for a self-contained exposition.
\par
Similarly, \cite[Theorem 1.1]{DGH} suggests that Johnson's original definition of amenability is too strong to deal with measure algebras. In \cite{Run}, the first-named author extended the
notion of Connes-amenability to the class of dual Banach algebras. This class includes --- besides $W^\ast$-algebras --- all measure algebras and all algebras ${\cal B}(E)$ for
a reflexive Banach space $E$. In \cite{Run2}, we proved that a locally compact group $G$ is amenable if and only if $M(G)$ is Connes-amenable.
\par
We shall now introduce a hybrid of operator amenability and Connes-amenability, which will turn out to be the ``right'' notion of amenability for the reduced Fourier--Stieltjes algebra
in the sense that it singles out precisely the amenable, locally compact groups.
\begin{definition}
A completely contractive Banach algebra is called {\it dual\/} if $\A^\ast$ has a closed $\A$-submodule $\A_\ast$ such that $\A = (\A_\ast)^\ast$.
\end{definition} 
\begin{remark}
In general, there is no need for $\A_\ast$ to be unique.
\end{remark}
\begin{examples}
\item If $\A$ is a dual Banach algebra in the sense of \cite[Definition 1.1]{Run}, then $\max \A$ is a dual, completely contractive Banach algebra.
\item Every $W^\ast$-algebra is a dual, completely contractive Banach algebra.
\item For any locally compact group $G$, the Fourier--Stieltjes algebras $B(G)$ and $B_r(G)$ are dual, completely contractive Banach algebras.
\end{examples}
\begin{definition}
Let $\A$ be a dual, completely contractive Banach algebra. A dual operator $\A$-bimodule $E$ is called {\it normal\/} if, for each $x \in E$, the maps
\[
  \A \to E, \quad a \mapsto \left\{ \begin{array}{c} a \cdot x \\ x \cdot a \end{array} \right.
\]
are $w^\ast$-continuous.
\end{definition}
\begin{sloppy}
\begin{definition}
A dual, completely contractive Banach algebra $\A$ is called {\it operator Connes-amenable\/} if every $w^\ast$-continuous, completely bounded derivation from
$\A$ into a normal, dual operator $\A$-bimodule is inner.
\end{definition}
\end{sloppy}      
\par
For the reduced Fourier--Stieltjes algebra, we obtain:
\begin{theorem} \label{opConnes}
The following are equivalent for a locally compact group $G$:
\begin{items}
\item $G$ is amenable.
\item $B_r(G)$ is operator Connes-amenable.
\end{items}
\end{theorem}
\begin{proof}
(i) $\Longrightarrow$ (ii): By \cite[Theorem 3.6]{Rua}, $A(G)$ is operator amenable. The $w^\ast$-density of $A(G)$ in $B_r(G)$ then yields the operator Connes-amenability
of $B_r(G)$ (compare \cite[Proposition 4.2(i)]{Run}).
\par
(ii) $\Longrightarrow$ (i): The same argument as in the proof of \cite[Proposition 4.1]{Run} yields that $B_r(G)$ has an identity. Since $B_r(G)$ is a closed ideal of $B(G)$
by \cite[(2.16) Proposition]{Eym}, it follows that $B_r(G) = B(G)$ and thus $C^\ast_r(G) = \cstar(G)$. By \cite[(4.21) Theorem]{Pat}, this is equivalent to $G$ being amenable.
\end{proof}
\par
It is, of course, tempting to conjecture that $B_r(G)$ in Theorem \ref{opConnes}(ii) can be replaced by $B(G)$. The implication (i) $\Longrightarrow$ (ii) then still holds because
$B(G) = B_r(G)$ for amenable $G$. The argument used to establish the converse, however, does no longer work for $B(G)$ instead of $B_r(G)$. As well shall now see, not only the proof no longer works, but
the statement becomes false: there are non-amenable, locally compact groups for which $B(G)$ is operator Connes-amenable.
\begin{lemma} \label{nicolem}
Let $G$ be a locally compact group. Then $A_{\cal F}(G)$ is operator amenable.
\end{lemma}
\begin{proof}
Since $aG$ is compact, we have $B(aG) = A(aG)$. Since $aG$ is amenable, $A(aG) = B(aG)$ is operator amenable by \cite[Theorem 3.6]{Rua}. By Proposition \ref{nicoprop}, the completely contractive Banach algebras $B(aG)$ and $A_{\cal F}(G)$ 
are completely isometrically isomorphic. Hence, $A_{\cal F}(G)$ is operator amenable.
\end{proof}
\begin{remark}
Should our conjecture that $B(G)$ is operator amenable only for compact $G$ be correct, then Lemma \ref{nicolem} would yield immediately that $A(G) \subsetneq A_{\PIF}(G)$ for non-compact, amenable $G$: otherwise,
we would have a short exact sequence 
\[
  \{ 0 \} \to A(G) \to B(G) \to A_{\cal F}(G) \to \{ 0 \}
\]
of completely contractive Banach algebras whose endpoints are operator amenable. The straightforward analogue of a hereditary property of amenability in the classical sense (\cite[Theorem 2.3.10]{LoA}) would then yield the
operator amenability of $B(G)$, which is impossible.
\end{remark}
\par
Recall that a $\cstar$-algebra $\A$ is called {\it residually finite-dimensional\/} if the family of finite-dimensional $^\ast$-representations of $\A$ separates the points of $\A$. For locally compact groups $G$,
the property of $\cstar(G)$ being residually finite-dimensional implies that $G$ is maximally almost periodic (\cite[Theorem 1.1]{SpW}), though the converse need not be true (\cite{Bek}).
\begin{theorem} \label{nicothm}
Let $G$ be a locally compact group such that $\cstar(G)$ is residually finite-dimensional. Then $B(G)$ is operator Connes-amenable.
\end{theorem}
\begin{proof}
By Lemma \ref{nicolem}, $A_{\cal F}(G)$ is operator amenable. Since $\cstar(G)$ is residually finite-dimensional, a simple Hahn--Banach argument shows that $A_{\cal F}(G)$ is $w^\ast$-dense in $B(G)$.
Then (the operator analogue of) \cite[Proposition 4.2(i)]{Run} yields the operator Connes-amenability of $B(G)$.
\end{proof}
\begin{example}
Let $\free_2$ denote the free group in two generators. Then $\cstar(\free_2)$ is residually finite-dimensional by \cite[Proposition VII.6.1]{Dav}, so that $B(\free_2)$ is operator Connes-amenable by
Theorem \ref{nicothm}. However, $\free_2$ is not amenable.
\end{example}
\dated
\vfill
\begin{tabbing}
{\it Second author's address\/}: \= Department of Mathematical and Statistical Sciences \kill
{\it First author's address\/}: \> Department of Mathematical and Statistical Sciences \\
\> University of Alberta \\
\> Edmonton, Alberta \\
\> Canada T6G 2G1 \\[\medskipamount]
{\it E-mail\/}: \> {\tt vrunde@ualberta.ca} \\[\medskipamount]
{\it URL\/}: \> {\tt http://www.math.ualberta.ca/$^\sim$runde/runde.html} \\[\bigskipamount]
{\it Second author's address\/}: \> Department of Pure Mathematics \\
\> University of Waterloo \\
\> Waterloo, Ontario \\
\> Canada N2L 3G1 \\[\medskipamount]
{\it E-mail\/}: \> {\tt nspronk@math.uwaterloo.ca} 
\end{tabbing} 

\begin{thebibliography}{00}
%
\begin{small}
%
\bibitem{Ars} {\sc G.\ Arsac}, Sur l'espace de Banach engendr\'e par les coefficients d'une repr\'esentations unitaire. {\it Publ.\ D\'ep.\ Math.\ (Lyon)\/} {\bf 13\/} (1976), 1--101. 
%
\bibitem{Bek} {\sc M.\ E.\ B.\ Bekka}, On the full $\cstar$-algebras of arithmetic groups and the congruence subgroup problem. {\it Forum Math.\/}\ {\bf 11\/} (1999), 705--715.
%
\bibitem{Bur} {\sc R.\ B.\ Burckel}, {\it Weakly Almost Periodic Functions on Semigroups\/}. Gordon and Breach, 1970.        
%
\bibitem{DGH} {\sc H.\ G.\ Dales}, {\sc F.\ Ghahramani}, and {\sc A. \ Ya.\ Helemski\u{\i}}, The amenability of measure algebras. {\it J.\ London Math.\ Soc.\/}\ (2) {\bf 66\/} (2002), 213--226.
%
\bibitem{Dav} {\sc K.\ R.\ Davidson}, {\it $\cstar$-Algebras by Example\/}. American Mathematical Society, 1996.
%
\bibitem{ER} {\sc E.\ G.\ Effros} and {\sc Z.-J.\ Ruan}, {\it Operator Spaces\/}. Oxford University Press, 2000. 
%
\bibitem{Eym} {\sc P.\ Eymard}, L'alg\`ebre de Fourier d'un groupe localement compact. {\it Bull.\ Soc.\ Math.\ France\/} {\bf 92\/} (1964), 181--236.
%
\bibitem{For} {\sc B.\ E.\ Forrest}, Amenability and weak amenability of the Fourier algebra. Preprint (2000).
%
\bibitem{Hel} {\sc A.\ Ya.\ Helemski\u{\i}}, Homological essence of amenability in the sense of A.\ Connes: the injectivity of the predual bimodule (translated from the Russion).
{\it Math.\ USSR--Sb.\/}\ {\bf 68\/} (1991), 555--566.
%
\bibitem{Joh1} {\sc B.\ E.\ Johnson}, Cohomology in Banach algebras. {\it Mem.\ Amer.\ Math.\ Soc.\/}\ {\bf 127\/} (1972). 
%
\bibitem{Joh2} {\sc B.\ E.\ Johnson}, Approximate diagonals and cohomology of certain annihilator Banach algebras. {\it Amer.\ J.\ Math.\/}\ {\bf 94\/} (1972), 685--698.   
%
\bibitem{Joh3} {\sc B.\ E.\ Johnson}, Non-amenability of the Fourier algebra of a compact group. {\it J.\ London Math.\ Soc.\/}\ (2)
{\bf 50\/} (1994), 361--374. 
%
\bibitem{JKR} {\sc B.\ E.\ Johnson}, {\sc R.\ V.\ Kadison}, and {\sc J.\ Ringrose}, Cohomology of operator algebras, III. {\it Bull.\ Soc.\ Math.\ France\/} 
{\bf 100\/} (1972), 73--79.
%
\bibitem{LL} {\sc A.\ T.-M.\ Lau} and {\sc R.\ J.\ Loy}, Amenability of convolution algebras. {\it Math.\ Scand.\/}\ {\bf 79\/} (1996), 283--296.      
%
\bibitem{LLW} {\sc A.\ T.-M.\ Lau}, {\sc R.\ J.\ Loy}, and {\sc G.\ A.\ Willis}, Amenability of Banach and $\cstar$-algebras on locally compact groups. 
{\it Studia Math.\/}\ {\bf 119\/} (1996), 161--178. 
%
\bibitem{Mia} {\sc T.\ Miao}, Decomposition of $B(G)$. {\it Trans.\ Amer.\ Math.\ Soc.\/}\ {\bf 581\/} (1999), 4675--4692. 
%
\bibitem{Pal} {\sc T.\ W.\ Palmer\/}, {\it Banach Algebras and the General Theory of $^\ast$-Algebras\/}, I. Cambridge University Press, 1994.           
%
\bibitem{Pat} {\sc A.\ L.\ T.\ Paterson}, {\it Amenability\/}. American Mathematical Society, 1988.
%
\bibitem{Pis} {\sc G.\ Pisier}, {\it An Introduction to the Theory of Operator Spaces\/}. Notes du Cours du Centre Emile Borel, 2000.     
%
\bibitem{Rua} {\sc Z.-J.\ Ruan}, The operator amenability of $A(G)$. {\it Amer.\ J.\ Math.\/}\ {\bf 117\/} (1995), 1449--1474.
%
\bibitem{LoA} {\sc V.\ Runde}, {\it Lectures on Amenability\/}. Lectures Notes in Mathematics {\bf 1774\/}, Springer Verlag, 2002. 
%
\bibitem{Run} {\sc V.\ Runde}, Amenability for dual Banach algebras. {\it Studia Math.\/}\ {\bf 148\/} (2001), 47--66. 
%
\bibitem{Run2} {\sc V.\ Runde}, Connes-amenability and normal, virtual diagonals for measure algebras, I. {\it J.\ London Math.\ Soc.\/}\ (to appear).
%
\bibitem{SpW} {\sc N.\ Spronk} and {\sc P.\ J.\ Wood}, Diagonal type conditions on group $\cstar$-algebras. {\it Proc.\ Amer.\ Math.\ Soc.\/}\ {\bf129\/} (2000), 609--616.
%
\bibitem{Tak} {\sc M.\ Takesaki}, {\it Theory of Operator Algebras\/}, I. Springer Verlag, 1979. 
%
\bibitem{Wal} {\sc M.\ E.\ Walter}, $\wstar$-algebras a nonabelian harmonic analysis. {\it J.\ Funct.\ Ana.\/}\ {\bf 11\/} (1972), 17--38.
%
\bibitem{Was} {\sc S.\ Wassermann}, On Tensor products of certain group $\cstar$-algebras. {\it J.\ Funct.\ Anal.\/}\ {\bf 23\/} (1976), 239--254.  
%
\bibitem{Wit} {\sc G.\ Wittstock} {\it et al.\/}, What are operator spaces? --- An online dictionary. URL: {\tt http://www.math.uni-sb.de/$^\sim$ag-wittstock/projekt2001.html}
(2001).   
%
\end{small} 
%
\end{thebibliography}
\end{document}